\title{Finite simple automorphism groups\\ of edge-transitive maps}

\author{Gareth A. Jones\\
School of Mathematical Sciences\\
University of Southampton\\
Southampton SO17  1BJ, UK\\
{\tt G.A.Jones@maths.soton.ac.uk}\\
}

\documentclass[12pt]{article}
\usepackage[latin1]{inputenc}
\usepackage{latexsym}
\usepackage{amsmath}
\usepackage{amssymb}
\usepackage{amsfonts}
\usepackage{color}
\usepackage{tikz}
\usepackage{enumerate}
\newtheorem{thm}{Theorem}[section]
\newtheorem{lemma}[thm]{Lemma}
\newtheorem{cor}[thm]{Corollary}
\newtheorem{prop}[thm]{Proposition}

\newcommand{\M}{\mathcal{M}}
\newcommand{\F}{\mathbb{F}}
\newcommand{\G}{\mathcal{G}}

\date{}
\begin{document} 

\newpage

\maketitle

\centerline{This paper is dedicated to the memory of Jan Saxl.}

\begin{abstract}
\noindent Building on earlier results for regular maps and for orientably regular chiral maps, we classify the non-abelian finite simple groups arising as automorphism groups of maps in each of the $14$ Graver--Watkins classes of edge-transitive maps.
\end{abstract}

\medskip


 \medskip
 
 
 \bigskip


\section{Introduction}

The principal result of this paper, Theorem~\ref{mainthmsimple}, classifies the non-abelian finite simple groups which can arise as the automorphism group of an edge-transitive map in each of the $14$ Graver--Watkins classes of such maps. Before stating the theorem precisely, we will discuss the background to this result in both finite group theory and topological graph theory.

In 1980 Mazurov asked in the Kourovka Notebook~\cite[Problem 7.30]{Kou} which finite simple groups have a generating set of three involutions, two of which commute (now often called a {\em Mazurov triple}). Since the classification of finite simple groups was announced around this time, his problem could be solved by case-by-case analysis of the alternating groups, the simple groups of Lie type, and the 26 sporadic simple groups. Nuzhin dealt with the alternating groups and those of Lie type in~\cite{Nuz90, Nuz92, Nuz97a, Nuz97}. In~\cite{Maz} Mazurov gave an elegant, unified and largely computer-free discussion of which sporadic simple groups have this property, together with a summary of how mathematicians such as Ershov and Nevmerzhitskaya, Nuzhin and Timofeenko, Abasheev, and Norton earlier dealt with various individual sporadic groups, mainly by using computers. The announced solution was that every non-abelian finite simple group has a Mazurov triple, with a small list of exceptions (including several infinite families of low Lie rank).

During the 1980s an algebraic approach to the theory of maps on surfaces was being developed, much of it inspired by Chapter 8 of the book~\cite{CM} by Coxeter and Moser and by Tutte's paper~\cite{Tut}. The most symmetric maps are the regular maps, those for which the automorphism group acts transitively on vertex-edge-face flags. A group arises as the automorphism group of such a map if and only if it has a Mazurov triple, so the solution of the Kourovka Notebook problem appeared also to give a classification of the finite simple automorphism groups of regular maps, and indeed the resulting classification was published in several papers: see \v Sir\'a\v n's excellent survey~\cite{Sir13}, for example. 

If it is known or suspected that a particular group $G$ can be realised as the automorphism group of a regular map, then it is natural to try to determine the number $n(G)$ of such maps, up to isomorphism. These maps correspond bijectively to the orbits of ${\rm Aut}\,G$ on Mazurov triples; this action is semiregular, so
\[n(G)=m(G)/|{\rm Aut}\,G|\]
where $m(G)$ is the number of Mazurov triples in $G$. The latter can be found, if $G$ is not too large, by hand or by a computer search, giving $n(G)$.

In 2016 Martin Ma\v caj~\cite{Maca} used GAP to apply this method to a number of finite simple groups $G$, and found that $n(G)=0$ for $G=U_4(3)$ and $U_5(2)$, even though these groups did not appear in the published lists of exceptions. His result was verified by Marston Conder and Matan Ziv-Av, using Magma and GAP. The earlier work on Mazurov triples in unitary groups was re-examined, and it became clear that a minor error, subsequently corrected in~\cite{Nuz19}, had led to the omission of these two groups. In consequence, we now have the following amended theorem (using ATLAS notation~\cite{ATLAS}):

\begin{thm}\label{Mthm}
A non-abelian finite simple group is generated by a Mazurov triple if and only if it is not isomorphic to a group in the following  set $\mathbb M$, where $q$ is any prime power and $e\ge 1$:
\[\{L_3(q),\, U_3(q), L\,_4(2^e),\, U_4(2^e),\, U_4(3),\, U_5(2),\, A_6,\, A_7,\, M_{11},\, M_{22},\, M_{23},\, McL\}.\]
\end{thm}

Note that the exceptions include the groups $L_2(7)\cong L_3(2)$, $L_2(9)\cong A_6$, $A_8\cong L_4(2)$, the symplectic group $S_4(3)\cong U_4(2)$ and the orthogonal group $O_6^-(3)\cong U_4(3)$. This theorem immediately implies the following:

\begin{cor}\label{Mcor}
A non-abelian finite simple group is isomorphic to the automorphism group of a regular map if and only if it not isomorphic to a group in the set $\mathbb M$. 
\end{cor}

By definition, a map $\M$ is edge-transitive if ${\rm Aut}\,\M$ acts transitively on its edges. In 1997 Graver and Watkins~\cite{GW} partitioned edge-transitive maps $\M$ into $14$ classes, distinguished by the isomorphism class of the quotient map ${\mathcal M}/{\rm Aut}\,{\mathcal M}$ (see the first column of Table~\ref{FSGpsRealised} for a list of the classes, and \S\ref{etrans} for an explanation of the classification); in that year, Wilson gave a similar classification in~\cite{Wil97}. These classes $T$ correspond bijectively to the $14$ isomorphism classes of maps ${\mathcal N}(T)\cong{\mathcal M}/{\rm Aut}\,{\mathcal M}$ with one edge; they include class $1$, consisting of the regular maps, and class $2^P{\rm ex}$, consisting of the orientably regular chiral maps (those non-regular orientable maps with an arc-transitive orientation-preserving automorphism group), together with others such as class~3, consisting of the just-edge-transitive maps (those for which ${\rm Aut}\,\M$ is transitive on edges but not on vertices or faces).

The groups realised as automorphism groups in class~1 are those with Mazurov triples, or equivalently those arising as quotients of the free product $\Gamma=V_4*C_2$. Similarly, one can show that those realised in the other edge-transitive classes $T$ are the quotients of certain finitely presented `parent groups' $N(T)$, with the additional requirement that they should not possess any `forbidden automorphisms', which would make the corresponding maps too symmetric to be members of $T$. This gives, for each of the $14$ classes $T$, necessary and sufficient conditions for a group to be a member of the set
\[\G(T):=\{G\mid G\cong{\rm Aut}\,\M\;\hbox{for some map}\;\M\in T\}\]
of groups realised in $T$. As a simple example, for the class $T=2^P{\rm ex}$ of orientably regular chiral maps we have
\[N(T)=\langle X, Y\mid Y^2=1\rangle\cong C_{\infty}*C_2,\]
so ${\mathcal G}(2^P{\rm ex})$ consists of those groups $G$ generated by elements $x$ (rotation around a vertex) and $y$ (rotation around the midpoint of an incident edge) with $y^2=1$; the forbidden automorphism in this case is the one inverting $x$ and fixing $y$ (equivalently, inverting both), since the presence of this automorphism (corresponding to reflection in the incident edge) would make a map regular, in class $1$ rather than $2^P{\rm ex}$.

One can therefore examine various classes of groups, such as the finite simple groups, to see which of them are realised as automorphism groups in each of the classes $T$. Here useful results on generators of non-abelian finite simple groups have been provided by Malle, Saxl and Weigel in~\cite{MSW}; for instance, they show that every such group except $U_3(3)$ is generated by three involutions, a result relevant to class $T=2$, with $N(T)\cong C_2*C_2*C_2$. This task is eased by the presence of a group $\Omega\cong S_3$ of operations on maps, introduced by Wilson in~\cite{Wil}, which permutes the 14 classes, preserving the sets $\G(T)$; it has six orbits on the classes (corresponding to the rows of Table~\ref{FSGpsRealised} and of Figure~\ref{basicmaps}), four of length $3$ and two of length $1$ (the classes $1$ and $3$), so it is sufficient to determine $\G(T)$ for just six classes, representatives of the six orbits of $\Omega$. 

A further simplification arises from the fact that certain epimorphisms $N(T)\to N(T')$ for $T'=1$ and $2^P{\rm ex}$ show that groups in $\G(T')$ are also in $\G(T)$ for various other classes $T$.
This increases the importance of these two classes $T'$, adding to their intrinsic and historical  importance in the theory of maps. We have already discussed $\G(1)$, the set of groups with Mazurov triples. There has also been significant interest in $\G(2^P{\rm ex})$, and in 2016 Leemans and Liebeck announced (at SODO 2016, Queenstown, NZ, in the more general context of abstract polytopes) a classification of the non-abelian finite simple groups in this set:

\begin{thm}\label{Lthm}
A non-abelian finite simple group is generated by an element  $x$ and an involution $y$, with no automorphism inverting $x$ and fixing $y$, if and only if it is not isomorphic to a group in the following set $\mathbb L$, where $q$ is any prime power:
\[\{L_2(q),\, L_3(q),\, U_3(q),\, A_7\}.\]
\end{thm}

As in the case of Theorem~\ref{Mthm}, there is an immediate corollary:

\begin{cor}\label{Lcor}
A non-abelian finite simple group is isomorphic to the automorphism group of an orientably regular chiral map if and only if it not isomorphic to a group in the set $\mathbb L$. 
\end{cor}

In~\cite{LL} Leemans and Liebeck proved these results in one direction, showing that all non-abelian finite simple groups except those in $\mathbb L$ have such a generating pair $x, y$, and are therefore automorphism groups of orientably regular chiral maps, that is, members of $\G(2^P{\rm ex})$. For the converse, it is easy to show that $L_2(q)$ and $A_7$ have no such generating pairs, and in~\cite{BC19} d'Azevedo Breda and Catalano proved this for the groups $L_3(q)$ and $U_3(q)$, thus completing the proof of Theorem~\ref{Lthm}.

\smallskip

In dealing with the classes $T\ne 1, 2^P{\rm ex}$ it is thus sufficient to restrict attention to $T=2, 3, 4$ and $5$ and to the groups in the following set ${\mathbb L}\cup{\mathbb M}$:
\[\{L_2(q),\, L_3(q),\, U_3(q), L\,_4(2^e),\, U_4(2^e),\, U_4(3),\, U_5(2),\]
\[A_6,\, A_7,\, M_{11},\, M_{22},\, M_{23},\, McL\}.\]
This is because the results for $T=1$ and $2^P{\rm ex}$, together with the simplifications mentioned above, imply that all other non-abelian finite simple groups are in $\G(T)$ for all classes $T$. The potential exceptions in ${\mathbb L}\cup{\mathbb M}$ can be treated individually, leading to the main result of this paper:

\begin{thm}\label{mainthmsimple}
A non-abelian finite simple group is isomorphic to the automorphism group of a map in an edge-transitive class $T$ if and only if it is not isomorphic to one of the groups listed in the corresponding row of Table~\ref{FSGpsRealised}.
\end{thm}

\begin{table}[ht]
\centering
\begin{tabular}{| p{2.6cm} | p{7.9cm}|}
\hline
Class $T$ & Non-abelian finite simple groups $G\not\in\G(T)$   \\
\hline\hline
$1$ & $L_3(q), U_3(q), L_4(2^e), U_4(2^e), U_4(3), U_5(2)$,\\
& $A_6, A_7, M_{11}, M_{22}, M_{23}, McL$  \\
\hline
$2, 2^*, 2^P$  & $U_3(3)$  \\
\hline
$2\,{\rm ex}, 2^*{\rm ex}, 2^P{\rm ex}$  &$L_2(q), L_3(q), U_3(q), A_7$  \\
\hline
$3$ & --  \\
\hline
$4, 4^*,4^P$ & --  \\
\hline
$5, 5^*, 5^P$ & $L_2(q)$  \\
\hline

\end{tabular}
\caption{Non-abelian finite simple groups not in sets $\G(T)$.}
\label{FSGpsRealised}
\end{table}

In order to avoid misunderstandings, we repeat here that the members of ${\mathbb L}\cup{\mathbb M}$ appearing as exceptions in various rows of Table~\ref{FSGpsRealised} include the groups $L_2(7)\cong L_3(2)$, $L_2(9)\cong A_6$, $A_8\cong L_4(2)$, the symplectic group $S_4(3)\cong U_4(2)$ and the orthogonal group $O_6^-(3)\cong U_4(3)$.

After some preliminary results and techniques concerning maps are described in the next section, the strategy for proving Theorem~\ref{mainthmsimple} is to consider the various groups in ${\mathbb L}\cup{\mathbb M}$, starting with the alternating groups, and to determine those sets $\G(T)$ to which they belong. 


\bigskip

\noindent{\bf Acknowledgements} The author is grateful to David Craven, Kay Magaard and Chris Parker for many useful hints about finite simple groups, to Dimitri Leemans for sharing details of his work with Martin Liebeck on automorphism groups of chiral polyhedra, to Marston Conder, Martin Ma\v caj and Matan Ziv-Av for their computations, and to Jan Saxl for fifty years of friendship and for his many valuable contributions to algebra.


\section{Algebraic map theory}\label{algthymaps}

\begin{figure}[h!]
\begin{center}
\begin{tikzpicture}[scale=0.5, inner sep=0.8mm]

\node (c) at (0,0) [shape=circle, fill=black] {};
\node (d) at (8,0) [shape=circle, fill=black] {};
\draw [thick] (c) to (d);
\draw [thick] (c) to (1,-3);
\draw [thick] (c) to (1,3);
\draw [thick] (d) to (7,-3);
\draw [thick] (d) to (7,3);
\draw [thick] (c) to (-2.5,2.5);
\draw [thick] (c) to (-2.5,-2.5);
\draw [thick] (d) to (10.5,2.5);
\draw [thick] (d) to (10.5,-2.5);

\draw (c) to (1,0.5);
\draw (c) to (1,-0.5);
\draw (1,0.5) to (1,-0.5);
\draw (d) to (7,0.5);
\draw (d) to (7,-0.5);
\draw (7,0.5) to (7,-0.5);
\draw (c) to (0.8,0.8);
\draw (0.3,1) to (0.8,0.8);

\node at (-0.8,0) {$v$};
\node at (4,-0.4) {$e$};
\node at (4,2) {$f$};

\node at (1.5,0.5) {$\phi$};
\node at (6.2,0.5) {$\phi r_0$};
\node at (1.4,1.3) {$\phi r_1$};
\node at (1.9,-0.5) {$\phi r_2$};
\node at (5.9,-0.5) {$\phi r_0r_2$};

\end{tikzpicture}

\end{center}
\caption{Action of generators $r_i$ of $B$ on a flag $\phi=(v,e,f)$.} 
\bigskip
\label{flags}
\end{figure}
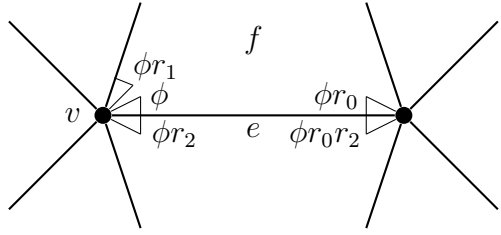

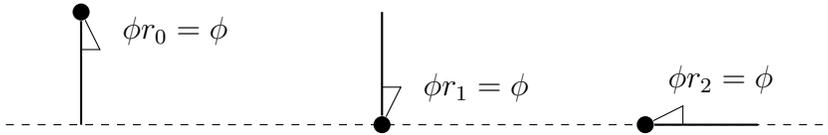
\begin{figure}[h!]
\begin{center}
\begin{tikzpicture}[scale=0.5, inner sep=0.8mm]

\draw [dashed] (-10,0) to (12,0);

\node (a) at (-8,3) [shape=circle, fill=black] {};
\draw [thick] (a) to (-8,0);
\draw (a) to (-7.5,2) to (-8,2);
\node at (-5.5,2.5) {$\phi r_0=\phi$};

\node (b) at (0,0) [shape=circle, fill=black] {};
\draw [thick] (b) to (0,3);
\draw (b) to (0.5,1) to (0,1);
\node at (2.5,1) {$\phi r_1=\phi$};

\node (c) at (7,0) [shape=circle, fill=black] {};
\draw [thick] (c) to (10,0);
\draw (c) to (8,0.5) to (8,0);
\node at (9,1.2) {$\phi r_2=\phi$};

\end{tikzpicture}

\end{center}
\caption{Flags fixed by generators $r_0, r_1$ and $r_2$.}
\label{fixedflags}
\end{figure}

This section outlines the algebraic theory of maps described in more detail elsewhere, such as in~\cite{JT}; see~\cite{Sir13} for regular maps, and~\cite{GT} for further background in topological graph theory. 
Any map $\mathcal M$ (possibly non-orientable or with non-empty boundary) determines a permutation representation of the group
\[\Gamma=\langle R_0, R_1, R_2\mid R_i^2=(R_0R_2)^2=1\rangle\cong V_4*C_2\]
on the set $\Phi$ of incident vertex-edge-face flags $\phi=(v,e,f)$ of $\mathcal M$. For any $\phi\in\Phi$ and $i=0, 1$ or $2$, there is at most one flag $\phi'\ne \phi$ with the same $j$-dimensional components as $\phi$ for $j\ne i$ (there may be none if $\phi$ is on the boundary of $\mathcal M$). Let $r_i$ be the permutation of $\Phi$ transposing each $\phi$ with $\phi'$ if the latter exists, and fixing $\phi$ otherwise, in which case we will call the incident edge a {\em free edge} or {\em semi-edge}. Figures~\ref{flags} and~\ref{fixedflags} illustrate these two cases, with the broken line in Figure~\ref{fixedflags} representing part of the boundary. Clearly $r_i^2=(r_0r_2)^2=1$, so there is a permutation representation
\[\theta:\Gamma\to {\rm Sym}\,\Phi\]
of $\Gamma$ on $\Phi$ given by $R_i\mapsto r_i$. The permutation group $B:=\langle r_0, r_1, r_2\rangle\le {\rm Sym}\,\Phi$ generated by $r_0, r_1$ and $r_2$ is called the {\em monodromy group\/} ${\rm Mon}\,{\mathcal M}$ of $\mathcal M$.

Conversely, any permutation representation of $\Gamma$ on a set $\Phi$ determines a map $\mathcal M$: the vertices, edges and faces are identified with the orbits on $\Phi$ of the subgroups $\langle R_1, R_2\rangle\cong D_{\infty}$,  $\langle R_0, R_2\rangle\cong V_4$ and  $\langle R_0, R_1\rangle\cong D_{\infty}$, with incidence given by non-empty intersection.

A map $\mathcal M$ is connected if and only if $\Gamma$ is transitive on $\Phi$, as we will always assume. The stabilisers in $\Gamma$ of flags $\phi\in\Phi$ then form a conjugacy class of subgroups $M\le\Gamma$, called {\em map subgroups}.

A map $\mathcal M$ is finite if and only if $M$ has finite index in $\Gamma$, and it has non-empty boundary if and only if some $r_i$ has a fixed point in $\Phi$, or equivalently some $R_i$ is contained in a conjugate of $M$. In particular, $\M$ is orientable and without boundary if and only if $M$ is contained in the even subgroup $\Gamma^+$ of index $2$ in $\Gamma$, consisting of the words of even length in the generators $R_i$.

The {\em automorphism group\/} $A={\rm Aut}\,{\mathcal M}$ of $\mathcal M$ is the centraliser in ${\rm Sym}\,\Phi$ of the monodromy group $B$. It is isomorphic to $N/M$ where $N:=N_{\Gamma}(M)$ is the normaliser  of $M$ in $\Gamma$ (see~\cite[Theorem~1]{JonADAM} or~\cite[Theorem~3.2]{PS}, for example). A map $\mathcal M$ is {\em regular\/} if $A$ is transitive on $\Phi$, which is equivalent to $B$ being a regular permutation group, and to $M$ being normal in $\Gamma$; in this case
$A\cong B\cong \Gamma/M$,
and $\Phi$ can be identified with $B$ so that the automorphism and monodromy groups are the left and right regular representations of this group. A map $\mathcal M$ is {\em edge-transitive\/} if $A$ is transitive on its edges; this is equivalent to $\Gamma=NE$ where $E:=\langle R_0, R_2\rangle\cong V_4$.

The (classical) {\em dual\/} $D({\mathcal M})$ of a map $\mathcal M$, with vertices and faces transposed, corresponds to the image of $M$ under the automorphism of $\Gamma$ fixing $R_1$ and transposing $R_0$ and $R_2$. The {\em Petrie dual\/} $P({\mathcal M})$ embeds the same graph as $\mathcal M$, but its faces are transposed with Petrie polygons, closed zig-zag paths turning alternately first right and first left at the vertices of $\mathcal M$; this corresponds to the automorphism of $\Gamma$ fixing $R_1$ and $R_2$ and transposing $R_0$ with $R_0R_2$. Both $D$ and $P$ preserve automorphism groups and regularity, but $D$ may change the embedded graph, and $P$ may change the underlying surface; in particular, if $\M$ is orientable, then $P(\M)$ is orientable if and only if the embedded graph is bipartite. The group $\Omega=\langle D, P\rangle\cong S_3$ of map operations, introduced by Wilson in~\cite{Wil}, permutes vertices, faces and Petrie polygons; it corresponds to the outer automorphism group ${\rm Out}\,\Gamma\cong {\rm Aut}\, E\cong S_3$ of $\Gamma$ acting on maps by permuting conjugacy classes of map subgroups~\cite{JT}.


\section{The edge-transitive classes}\label{etrans}

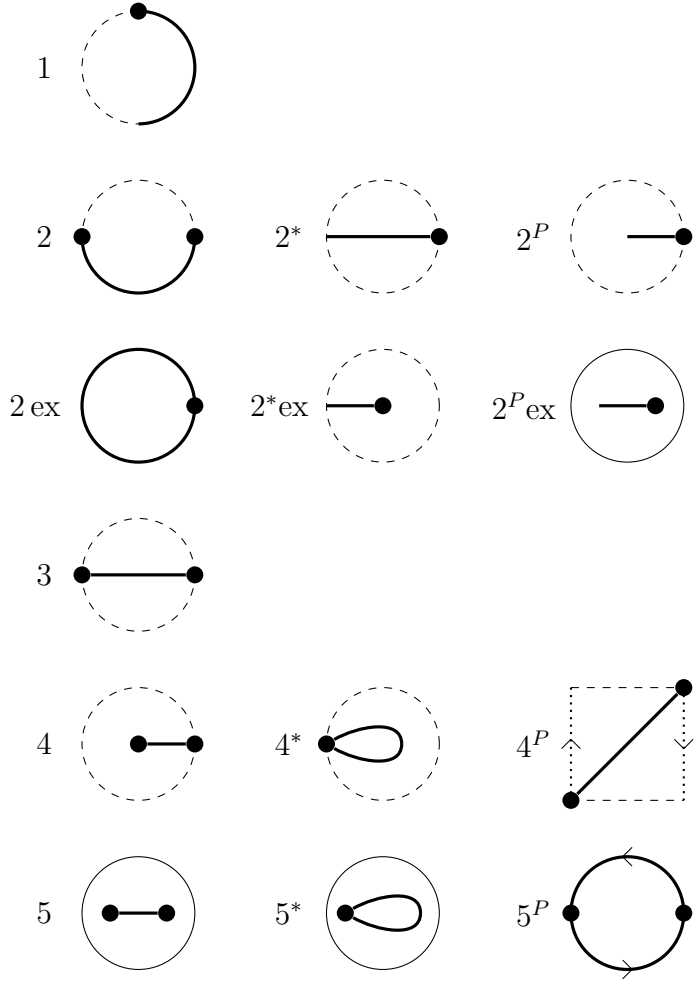
\begin{figure}[h!]
\begin{center}
\begin{tikzpicture}[scale=0.25, inner sep=0.8mm]

\node (a) at (0,48) [shape=circle, fill=black] {};
\draw [very thick] (0,42) arc (-90:90:3);
\draw [dashed] (0,48) arc (90:270:3);
\node at (-5,45) {$1$};


\node (b) at (-3,36) [shape=circle, fill=black] {};
\node (c) at (3,36) [shape=circle, fill=black] {};
\draw [very thick] (-3,36) arc (-180:0:3);
\draw [dashed] (3,36) arc (0:180:3);
\node at (-5,36) {$2$};


\node (d) at (16,36) [shape=circle, fill=black] {};
\draw [dashed] (16,36) arc (0:360:3);
\draw [very thick] (d) to (10,36);
\node at (8,36) {$2^*$};


\node (e) at (29,36) [shape=circle, fill=black] {};
\draw [dashed] (29,36) arc (0:360:3);
\draw [very thick] (e) to (26,36);
\node at (21,36) {$2^P$};


\node (f) at (3,27) [shape=circle, fill=black] {};
\draw [very thick] (3,27) arc (0:360:3);
\node at (-5.5,27) {$2\,{\rm ex}$};


\node (g) at (13,27) [shape=circle, fill=black] {};
\draw [dashed] (16,27) arc (0:360:3);
\draw [very thick] (g) to (10,27);
\node at (7.5,27) {$2^*$ex};


\node (h) at (27.5,27) [shape=circle, fill=black] {};
\draw [thin] (29,27) arc (0:360:3);
\draw [very thick] (h) to (24.5,27);
\node at (20.5,27) {$2^P$ex};


\node (i) at (-3,18) [shape=circle, fill=black] {};
\node (j) at (3,18) [shape=circle, fill=black] {};
\draw [very thick] (i) to (j);
\draw [dashed] (3,18) arc (0:360:3);
\node at (-5,18) {$3$};


\node (k) at (0,9) [shape=circle, fill=black] {};
\node (l) at (3,9) [shape=circle, fill=black] {};
\draw [very thick] (k) to (l);
\draw [dashed] (3,9) arc (0:360:3);
\node at (-5,9) {$4$};


\draw [dashed] (16,9) arc (0:360:3);
\node (m) at (10,9) [shape=circle, fill=black] {};
\draw [very thick] (m) to [out=30,in=90] (14,9);
\draw [very thick] (m) to [out=-30,in=-90] (14,9);
\node at (8,9) {$4^*$};


\draw [dashed] (23,12) to (29,12);
\draw [dashed] (23,6) to (29,6);
\draw [thick] [dotted] (23,12) to (23,6);
\draw [thick] [dotted] (29,12) to (29,6);
\node (n) at (23,6) [shape=circle, fill=black] {};
\node (o) at (29,12) [shape=circle, fill=black] {};
\draw [very thick] (n) to (o);
\draw [thin] (22.5,8.75) to (23,9.25);
\draw [thin] (23.5,8.75) to (23,9.25);
\draw [thin] (28.5,9.25) to (29,8.75);
\draw [thin] (29.5,9.25) to (29,8.75);
\node at (21,9) {$4^P$};


\node (a) at (-1.5,0) [shape=circle, fill=black] {};
\node (b) at (1.5,0) [shape=circle, fill=black] {};
\draw [thin] (3,0) arc (0:360:3);
\draw [very thick] (a) to (b);
\node at (-5,0) {$5$};


\draw [thin] (16,0) arc (0:360:3);
\node (c) at (11,0) [shape=circle, fill=black] {};
\draw [very thick] (c) to [out=30,in=90] (15,0);
\draw [very thick] (c) to [out=-30,in=-90] (15,0);
\node at (8,0) {$5^*$};


\node (d) at (23,0) [shape=circle, fill=black] {};
\node (e) at (29,0) [shape=circle, fill=black] {};
\draw [very thick] (29,0) arc (0:360:3);
\draw (25.75,3) to (26.25,3.5);
\draw (25.75,3) to (26.25,2.5);
\draw (26.25,-3) to (25.75,-3.5);
\draw (26.25,-3) to (25.75,-2.5);
\node at (21,0) {$5^P$};

\end{tikzpicture}

\end{center}
\caption{The basic maps $\mathcal N(T)$ for the $14$ edge-transitive classes $T$}
\label{basicmaps}
\end{figure}

For any map $\M$, the quotient map $\M/{\rm Aut}\,\M$ can be found by taking a fundamental region $\mathcal F$ for ${\rm Aut}\,\M$ and identifying equivalent boundary points. When $\M$ is edge-transitive,  $\M/{\rm Aut}\,\M$ must be isomorphic to one of the 14 basic maps ${\mathcal N}(T)$ with a single edge, enabling one to determine the class $T$ to which $\M$ belongs. These basic maps are shown in Figure~\ref{basicmaps}, each labelled with the Graver--Watkins symbol for the class $T$.

A map $\M$ is edge-transitive if and only if $N_{\Gamma}(M)$ acts transitively on the cosets of $E=\langle R_0, R_2\rangle$ in $\Gamma$; this is equivalent to $\Gamma=N_{\Gamma}(M)E$, that is, $E$ acting transitively on the cosets of $N_{\Gamma}(M)$ in $\Gamma$. This condition implies that $|\Gamma:N_{\Gamma}(M)|\le|E|=4$. A simple analysis of involutions in $S_4$ shows that there are, up to equivalence, $14$ permutation representations of $\Gamma$ in which $E$ acts transitively. Constructing the maps $\mathcal N$ corresponding to these representations, we have the following:

\begin{thm}
There are $14$ classes $T$ of edge-transitive maps, corresponding to the $14$ conjugacy classes of subgroups $N$ of $\Gamma$ satisfying $\Gamma=NE$. These are the map subgroups for the maps ${\mathcal N}(T)$ shown in Figure~\ref{basicmaps}. \hfill$\square$
\end{thm}

The dualities $D$ and $P$ leave the classes $T=1$ and $3$ invariant, whereas in each of the other four rows of Figure~\ref{basicmaps}, $D$ transposes the first and second classes while fixing the third, while $P$ transposes the second and third while fixing the first. Thus each of the six rows is an orbit of $\Omega$. We will write $2^{\sigma}$ to denote an arbitrary class in the orbit containing $2$, where $\sigma$ denotes $\emptyset$ (the empty symbol), $*$ or $P$, and similarly for the orbits containing the classes $2\,{\rm ex},$ $4$ and $5$. Let
\[\G(T):=\{G\cong {\rm Aut}\,\M\mid \M\in T\}\]
denote the set of groups realised as automorphism groups of maps in $T$. Since $\Omega$ permutes these sets in the same way as it permutes the classes $T$, in determining the groups in the various sets $\G(T)$ it is sufficient to consider one representative $T$ of each orbit of $\Omega$, such as  $T=1$, $2$, $2\,{\rm ex}$, $3$, $4$ and $5$.

The Reidemeister--Schreier process gives the following presentations for representatives of the orbits of $\Omega$ on the $14$ parent groups $N(T)$:

\begin{prop}\label{parents}
$N(1)=\Gamma=\langle R_0, R_1, R_2\mid R_i^2=(R_0R_2)^2=1\rangle$,
\[N(2)=\langle S_1=R_1, S_2=R_1^{R_0}, S_3=R_2\mid S_1^2=S_2^2=S_3^2=1\rangle,\]
\[N(2\,{\rm ex})=\langle S_1=R_2, S=R_0R_1\mid S_1^2=1\rangle,\]
\[N(3)=\langle S_0=R_1,\, S_1=R_1^{R_0},\, S_2=R_1^{R_2},\, S_3=R_1^{R_0R_2}\mid S_i^2=1\rangle,\]
\[N(4)=\langle S_1=R_1, S_2=R_1^{R_2}, S=(R_1R_2)^{R_0}\mid S_1^2=S_2^2=1\rangle,\]
\[N(5)=\langle S=R_1R_2, S'=S^{R_0}\mid - \rangle.\]
\end{prop}

Applying elements of $\Omega$, permuting $R_0, R_2$ and $R_0R_2$, gives presentations for the other parent groups. It is often more convenient to take
\begin{equation}
N(2^P{\rm ex})=\Gamma^+=\langle X=R_1R_2, Y=R_0R_2 \mid Y^2=1\rangle
\end{equation}
to represent the third orbit, since there is considerable information in the literature on orientably regular maps, the members of classes $1$ and $2^*{\rm ex}$.

If $T=1$ the quotients $G$ of $N(T)=\Gamma$ are the groups generated by elements $r_0, r_1$ and $r_2$ of order dividing $2$ (images of $R_0, R_1$ and $R_2$), with $r_0$ and $r_2$ commuting. 

If $T=2^P{\rm ex}$ then $N(T)=\langle X, Y\mid Y^2=1\rangle$, and the quotients are the groups $G$ generated by elements $x$ and $y$ with $y^2=1$. The corresponding map $\mathcal M$ is in class $T=2^P{\rm ex}$ (rather than class~$1$) if and only if $M$ is not normalised by the element $R_2\in\Gamma\setminus N(2^P{\rm ex})$. Since $R_2$, acting by conjugation, inverts $X$ and fixes $Y$, this is equivalent to $G$ not having an automorphism inverting $x$ and fixing $y$.  A similar restriction applies when $T=2\,{\rm ex}$ or $2^*{\rm ex}$.

The arguments in the other cases are similar. To summarise the results (see also~\cite[Condition~3.2]{STW}), the $14$ classes $T$ are listed in Table~\ref{forbidden}, with the second column giving the isomorphism type of each parent group $N(T)$, and the third column indicating any forbidden automorphisms; these are given in terms of their effect on generators $s_i, s, s'$ of quotients $G$, the images of generators $S_i, S, S'$ of cyclic free factors of $N(T)$ defined in Proposition~\ref{parents}.

\begin{table}[htb]
\centering
\begin{tabular}{|c|c|c|}
\hline
Class $T$&$N(T)$&forbidden automorphisms\\
\hline\hline
$1$&$V_4*C_2$&none\\
\hline
$2^{\sigma}$&$C_2*C_2*C_2$&$s_1$ and $s_2$ transposed, $s_3$ fixed\\
\hline
$2^{\sigma}{\rm ex}$&$C_2*C_{\infty}$&$s_1$ (or $y$)  fixed, $s$ (or $x$) inverted\\
\hline
$3$&$C_2*C_2*C_2*C_2$&double transpositions of generators $s_i$\\
\hline
$4^{\sigma}$&$C_2*C_2*C_{\infty}$&$s_1$ and $s_2$ transposed, $s$ inverted\\
\hline
$5^{\sigma}$&$C_{\infty}*C_{\infty}=F_2$&$s$ and $s'$ inverted, transposed or both\\ \hline
\end{tabular}
\caption{Parent groups and forbidden automorphisms}
\label{forbidden}
\end{table}

The following lemma will prove useful later in realising various groups in certain classes, by showing that in many cases it is sufficient to consider just two classes, namely $T=1$ and $2^P{\rm ex}$.

\begin{lemma}\label{regmapslemma}
{\rm(a)} If $G$ is a non-abelian group in $\G(1)$, then $G\in\G(T)$ for each class $T=2^{\sigma}$, $3$ or $4^{\sigma}$.
\vskip2pt
{\rm(b)} If $G$ is any group in $\G(2^P{\rm ex})$, then $G\in\G(T)$ for each class $T=2^{\sigma}{\rm ex}$, $4^{\sigma}$ or $5^{\sigma}$.
\end{lemma}

\noindent{\sl Proof.} (a) By the hypothesis, $G$ has generators $r_0, r_1$ and $r_2$ satisfying $r_i^2=(r_0r_2)^2=1$. Since $G$ is non-abelian, applying duality if necessary we may assume that $r_1r_2$ has order  $n>2$. It is sufficient to prove the result for $T=2, 3$ and $4$.

First let $T=2$. The presentations of $\Gamma$ and $N(2)$ show that there is an epimorphism $N(2)\to\Gamma$ given by $S_i\mapsto R_{i-1}$ for $i=1, 2, 3$. Composing this with the epimorphism $\Gamma\to G, R_i\mapsto r_i$ gives an epimorphism $\theta: N(2)\to G$ with $S_i\mapsto s_i:=r_{i-1}$ for $i=1, 2, 3$. Any automorphism of $G$ transposing $s_1$ and $s_2$, and fixing $s_3$, would send $s_1s_3=r_0r_2$, which has order dividing $2$, to $s_2s_3=r_1r_2$, which has order $n>2$. This is impossible, so $G\in\G(2)$. 

A similar argument applies if $T=4$, with $S_i\mapsto r_{i-1}$ as before for $i=1, 2$, and $S\mapsto s:=r_2$. If $T=3$ take $S_i\mapsto s_i:=r_1, r_2, r_0, r_2$ for $i=0,\ldots, 3$. Any automorphism of $G$ inducing a double transposition on $s_0,\ldots, s_3$ would imply that $(r_0r_1)^2=1$, so $G\in\G(3)$.

{\rm(b)} We use the same method as in (a), except that we now have an epimorphism $N(2^P{\rm ex})=\Gamma^+\to G,\, X\mapsto x, Y\mapsto y$, where $G$ has no automorphism inverting $x$ and fixing $y$ (so $x$ has order $n>2$). Composing this with the epimorphism $N(5)\to N(2^P{\rm ex}),\, S\mapsto X, S'\mapsto Y$ gives an epimorphism $\theta: N(5)\to A,\, S\mapsto x, S'\mapsto y$. An automorphism of $G$ cannot invert $x$ and $y$, or transpose $y$ with $x^{\pm 1}$, so $G\in\G(5)$. A similar argument applies when $T=4$, with $S\mapsto x$ and $S_1, S_2\mapsto y$.
\hfill$\square$

\medskip

Part (c) of the following lemma shows that class~$2$ can also play a similar role to that played by classes~$1$ and $2^P{\rm ex}$ in Lemma~\ref{regmapslemma}. An element $x$ of a group $G$ is {\em strongly real\/} if $x^a=x^{-1}$ for some involution $a\in G$, or equivalently, if $x$ is a product of at most two involutions in $G$.

\begin{lemma}\label{class2lemma}
{\rm (a)} If a group $G$ is generated by involutions $a, b$ and $c$, where $ab$, $ac$ and $bc$ do not all have the same order, then $G\in\G(2)$.
\vskip2pt
{\rm(b)} If a group $G$ is generated by elements $a, b$ and $c$ satisfying $abc=1$, where $a$ is an involution, $b$ is a product of two involutions, and no automorphism of $G$ inverts $c$, then $G\in\G(2)$.
\vskip2pt
{\rm(c)} If $G\in\G(2)$ then $G\in\G(T)$ for each class $T=2^{\sigma}$, $3$ or $4^{\sigma}$.
\end{lemma}

\noindent{\sl Proof.} (a) Without loss of generality we may assume that $ac$ and $bc$ have different orders. There is an epimorphism $N(2)\to G$ given by $S_i\mapsto s_i:=a, b$ or $c$ for $i=1, 2$ or $3$. Then no automorphism of $G$ can transpose $s_1$ and $s_2$ while fixing $s_3$, so $G\in\G(2)$.

(b) We can write $b=s_1s_2$ for involutions $s_1, s_2\in G$, and define $s_3=a$, giving involutions $s_1, s_2, s_3$ generating $G$. If an automorphism transposes $s_1$ and $s_2$ while fixing $s_3$, it inverts $a$ and $b$, so composing it with conjugation by $a$ gives an automorphism inverting $c$. However, no such automorphism exists, so $G\in\G(2)$.

(c) It is sufficient to prove (c) for $T=3$ and $4$. The group $N(4)$ can be obtained from $N(2)$ by taking $S=S_3$ and omitting the relation $S_3^2=1$, giving an epimorphism $N(4)\to N(2)$. Thus any quotient $G=\langle s_1, s_2, s_3\rangle$ of $N(2)$ is also a quotient of $N(4)$, with $s:=s_3$; the forbidden automorphisms are the same in each case, since $s$ has order dividing $2$, so $\G(2)\subseteq\G(4)$. For $T=3$ one can extend an epimorphism $N(2)\to G$ to $N(3)\to G$ by mapping $S_0$ to $s_3$; again, this introduces no further forbidden automorphisms.  \hfill$\square$

\begin{lemma}\label{class3lemma}
Suppose that a group $G$ is generated by involutions $s_0, s_1, s_2, s_3$, and that for at least two of the three partitions $ij\mid kl$ of $\{0, 1, 2, 3\}$ the products $s_is_j$ and $s_ks_l$ have different orders. Then $G\in\G(3)$.  \hfill$\square$
\end{lemma}

\begin{lemma}
Let $G=\langle x, y\mid y^2=1, \ldots\rangle\in\G(2^P{\rm ex})$, where the image $x$ of $X$ is strongly real. Then $G\in\G(T)$ for each edge-transitive class $T\ne 1$.
\end{lemma}

\noindent{\sl Proof.} An involution $a\in G$ inverts $x$, so $G=\langle s_1, s_2, s_3\rangle$ where $s_1=a$, $s_2=ax$ and $s_3=y$ all satisfy $s_i^2=1$, and hence $G$ is a quotient of $N(2)$. If an automorphism of $G$ transposes $s_1$ and $s_2$ and fixes $s_3$, it inverts $x$ and fixes $y$, contradicting chirality of the corresponding map. Thus no such automorphism exists, so $G\in\G(2)$. Then $G\in\G(T)$ for each $T=2^{\alpha}$, $3$ or $4^{\alpha}$ by Lemma~\ref{class2lemma}(c), and also for $T=2^{\alpha}{\rm ex}$ or $5^{\alpha}$ by Lemma~\ref{regmapslemma}(b), so $G\in\G(T)$ for each $T\ne 1$. \hfill$\square$

\medskip

As noted by Singerman in~\cite{Sin74}, the next result follows from work of Macbeath~\cite{Macb69}:

\begin{lemma}[Macbeath, Singerman]\label{L2qlemma}
If two elements generate $L_2(q)$ then some automorphism of that group inverts them both.
\end{lemma}

To deal with some groups we will use information and notation concerning their conjugacy classes, characters, maximal subgroups, automorphisms, etc found in the ATLAS~\cite{ATLAS}, together with the following formula, due to Frobenius~\cite{Fro}:

\begin{prop}\label{frobchi}
Let $\mathcal A$, $\mathcal B$ and $\mathcal C$ be conjugacy classes in a finite group $G$. Then the number of solutions of the equation $abc=1$ in $G$, with $a\in{\mathcal A}$, $b\in{\mathcal B}$ and $c\in{\mathcal C}$, is 
\begin{equation}\label{trianglechi}
\frac{|{\mathcal A}|\cdot|{\mathcal B}|\cdot|{\mathcal C}|}{|G|}\sum_{\chi}
\frac{\chi(a)\chi(b)\chi(c)}{\chi(1)}
\end{equation}
where $\chi$ ranges over the irreducible complex characters of $G$. \hfill$\square$
\end{prop}

For arguments involving Singer cycles and their eigenvalues, we will need the following result:

\begin{lemma}\label{priminv}
Let $\lambda$ be a primitive root in a finite field $\F$. If $\lambda$ is conjugate to $\lambda^{-1}$ under ${\rm Gal}\,\F$ then $|\F|\le 4$.
\end{lemma}

\noindent{\sl Proof.} Let $|\F|=p^e$ where $p$ is prime, so that  ${\rm Gal}\,\F$ is generated by the Frobenius automorphism $t\mapsto t^p$, which has order $e$. If $\lambda$ is conjugate to its inverse, then $\lambda^{p^f}=\lambda^{-1}$ and hence $p^f+1$ is divisible by the multiplicative order $p^e-1$ of $\lambda$, for some $f=0, 1, \ldots, e-1$. Hence $p^e-1\le p^{e-1}+1$, so $p^{e-1}(p-1)\le 2$, and the result follows immediately.  \hfill$\square$


\section{Alternating groups}

Here we will use a classic theorem of Jordan, stated as Theorem~13.9 in~\cite{Wie}. It is not stated explicitly by Jordan, but it follows immediately from his results in \cite{Jor1, Jor2} (see~\cite{Jon14} for historical details and a generalisation, removing the primality condition).

\begin{lemma}\label{Jordan}
If $G$ is a primitive subgroup of $S_n$, containing a cycle of prime length with at least three fixed points, then $G\ge A_n$.
\end{lemma}

\begin{thm}\label{altgps}
Let $n\ge 5$ and let $T$ be an edge-transitive class of maps. Then $A_n\in\G(T)$ if and only if one of the following holds:
\begin{itemize}
\item $T=1$, and $n=5$ or $n\ge 9$;
\item $T=2^{\sigma}$, $3$ or $4^{\sigma}$ for some $\sigma$;
\item $T=2^{\sigma}{\rm ex}$ for some $\sigma$, and $n\ge 8$;
\item $T=5^{\sigma}$ for some $\sigma$, and $n\ge 7$.
\end{itemize}
\end{thm}

\noindent{\sl Proof.} Recall that, in order to determine which sets $\G(T)$ contain a specific group $G$, it is sufficient to consider the classes $T=1$, $2^*{\rm ex}$, $2$, $3$, $4$ and $5$. The only alternating groups in the set ${\mathbb L}\cup{\mathbb M}$ are $A_5\;(\cong L_2(5))$, $A_6\;(\cong L_2(9))$, $A_7$ and $A_8\;(\cong L_4(2))$, so $A_n\in\G(T)$ for all $T$ if $n>8$. We now consider the values $n=5,\ldots, 8$.

Since $A_5\in\G(1)$ (as the automorphism group of the regular embedding of $K_6$ in the projective plane, for example), it follows from Lemma~\ref{regmapslemma}(a) that $A_5\in\G(T)$ for all $T=2^{\sigma}, 3$ and $4^{\sigma}$.  On the other hand, Lemma~\ref{L2qlemma} shows that $A_5\not\in\G(T)$ if $T=2^{\sigma}{\rm ex}$ or $5^{\sigma}$, and the same applies to $A_6$.

We now consider the values $n=6$, $7$ and $8$ for the various classes. Firstly, $A_n\not\in\G(1)$ for such $n$ by Corollary~\ref{Mcor}.

If $T=2$ then $N(T)=\langle S_1, S_2,S_3\mid S_i^2=1\rangle\cong C_2*C_2*C_2$. In $A_6$ take $s_1=(1,2)(3,4)$, $s_2=(2,6)(4,5)$ and $s_3=(2,3)(4,5)$, so that
\[s_1s_2=(1,6,2)(3,5,4),\quad s_1s_3=(1,3,5,4,2)\quad{\rm and}\quad s_2s_3=(2,6,3).\]
In $A_7$, take $s_1=(1,2)(3,4)$, $s_2=(2,6)(5,7)$ and $s_3=(2,3)(4,5)$, so that
\[s_1s_2=(1,6,2)(3,4)(5,7),\quad s_1s_3=(1,3,5,4,2)\;\;{\rm and}\;\; s_2s_3=(2,6,3)(4,5,7).\]
In $A_8$ take $s_1=(1,2)(3,4)(5,6)(7,8)$, $s_2=(1,3)(4,6)$ and $s_3=(3,4)(6,7)$, so that
\[s_1s_2=(1,2,3,6,5,4)(7,8),\quad s_1s_3=(1,2)(5,7,8,6)\;\;{\rm and}\;\; s_2s_3=(1,4,7,6,3).\]
In each case, $G:=\langle s_1, s_2, s_3\rangle$ is a primitive subgroup of $A_n$: if $n=7$ then $G$ is transitive of prime degree, and if $n=6$ or $8$ it is doubly transitive, in the latter case because $\langle(s_1s_2)^2, s_3\rangle$ fixes $8$ and is transitive on the remaining points. Jordan's Theorem (Lemma~\ref{Jordan}), applied to a suitable $3$- or $5$-cycle, then gives $G=A_n$. In each case $s_1s_3$ and $s_2s_3$ have different orders, so $G$ has no automorphism transposing $s_1$ and $s_2$ and inverting (equivalently fixing) $s_3$. Thus $A_6, A_7, A_8\in\G(2)$, and the same follows for $T=2^{\sigma}, 3$ and $4^{\sigma}$ by Lemma~\ref{class2lemma}(c).

If $T=2^{\sigma}{\rm ex}$ and $5^{\sigma}$ then $A_6\not\in\G(T)$ by Lemma~\ref{L2qlemma} and the isomorphism with $L_2(9)$.  Corollary~\ref{Lcor} gives $A_7\not\in\G(2^{\sigma}{\rm ex})$, whereas $A_7\in\G(5^{\sigma})$ since the pair $(1,2,3,4,5)$ and $(1,6,7)(2,4,5)$ generate $A_7$ and are not inverted in ${\rm Aut}\,A_7=S_7$. Corollary~\ref{Lcor} also gives $A_8\in\G(2^*{\rm ex})$, so $A_8\in\G(T)$ for $T=2^{\sigma}{\rm ex}$ and $5^{\sigma}$ by Lemma~\ref{regmapslemma}(b) \hfill$\square$

\medskip

In fact, one can extend this argument to show that for each class $T$ the perfect double cover $2.A_n$ of $A_n$ is in $\G(T)$ for all sufficiently large $n$ ($n\ge 9$ is enough when $T\ne 1$); see~\cite[\S10]{JonArXiv} for details. Similar arguments in~\cite[\S7]{JonArXiv} show that $S_n\in\G(T)$ for all $T$ provided $n\ge6$.


\section{Sporadic simple groups}

Here we deal with the sporadic simple groups in ${\mathbb L}\cup{\mathbb M}$:

\begin{thm}
The Mathieu groups $M_{11}$, $M_{22}$, $M_{23}$ and the McLaughlin group $McL$ are in $\G(T)$ for all edge-transitive classes $T\ne 1$, but they are not in $\G(1)$.
\end{thm}

\noindent{\sl Proof.} Corollary~\ref{Mcor} shows that these groups are not in $\G(1)$. Corollary~\ref{Lcor} shows that they are in $\G(2^P{\rm ex})$, so by Lemma~\ref{regmapslemma}(b) they are in $\G(T)$ for all classes $T=2^{\sigma}{\rm ex}$, $4^{\sigma}$ and $5^{\sigma}$.

For $T=2^{\sigma}$ or $3$ we use Lemma~\ref{class2lemma}(b) and (c). The group $G=M_{11}$ has unique conjugacy classes $2A$ and $4A$ of elements of order $2$ and $4$, and two mutually inverse classes $11A$ and $11B$ of elements of order $11$. Proposition~\ref{frobchi} and the character table of $G$ in~\cite{ATLAS} show that there are elements $a, b, c\in G$, in classes $2A$, $4A$ and $11A$, such that $abc=1$. The maximal subgroups of $G$ are listed in~\cite{ATLAS}, and none contains elements of orders $4$ and $11$, so $\langle a, b\rangle=G$. A second application of Proposition~\ref{frobchi}, or equivalently the existence of subgroups isomorphic to $D_4$, shows that there are involutions $s_1, s_2\in G$ with $s_1s_2=b$. Since $G$ has no outer automorphisms, no automorphism can invert $c$, so Lemma~\ref{class2lemma} shows that $G\in\G(T)$ for all classes $T=2^{\sigma}$ and $3$.

Essentially the same argument gives the result for $M_{23}$, but now with $a, b$ and $c$ in classes $2A$, $8A$ and $23A$, and also for $McL$, using classes $2A$, $12A$ and $11A$. (Even though $McL$ has outer automorphisms, they do not transpose the mutually inverse classes $11A$ and $11B$.)

A similar argument can be applied to $G=M_{22}$, but in this case the details are less straightforward. This group has one class $2A$ of involutions, one class $6A$ of elements of order $6$, and two mutually inverse classes $7A$ and $7B$ of elements of order $7$; the latter are not transposed in ${\rm Aut}\,G$, which contains $G$ with index $2$. Proposition~\ref{frobchi} shows that $G$ contains $12|G|$ triples $(a, b, c)$ of elements of orders $2$, $6$ and $7$ with $abc=1$.  By inspection of the list of maximal subgroups of $G$, each of these triples can generate one of the $2|G:H|$ subgroups $H\cong A_7$, or one of the $|G:H|$ subgroups $H\cong AGL_3(2)$, or $G$ itself. The uniqueness of the corresponding permutation diagram for $a$ and $b$ shows that $A_7$ is generated by $|S_7|=2|A_7|$ such triples, giving a total of $4|G|$ triples $(a, b, c)$ generating subgroups $H\cong A_7$. The affine group $AGL_3(2)$ is a semidirect product $V_8\rtimes GL_3(2)$, and triples $(a, b, c)$ generating this group map onto triples $(\overline a, \overline b, \overline c)$ of type $(2, 3, 7)$ generating $GL_3(2)=L_3(2)$. There are $2|L_3(2)|$ such triples, each lifting to $16$ triples $(a, b, c)$, giving a total of $|G:H|.16.2|L_3(2)|=4|G|$ triples $(a, b, c)$ generating subgroups $H\cong AGL_3(2)$. This leaves $4|G|$ triples $(a, b, c)$ generating $G$ (and corresponding to a unique chiral pair of orientably regular maps of type $\{6, 7\}$ with automorphism group $G$). The element $b\in 6A$ in such a triple is strongly real (since $A_7$, and hence $G$, contains subgroups isomorphic to $D_6$), and the outer automorphisms of $G$ do not transposes the classes $7A$ and $7B$ containing $c$ and $c^{-1}$, so Lemma~\ref{class2lemma} applies as before. \hfill$\square$


\section{$L_2(q)$}

\begin{thm}
A simple group $L_2(q)$ is a member of $\G(T)$ if and only if either $T=1$ and $q\ne 7, 9$, or $T=2^{\sigma}, 3$ or $4^{\sigma}$ for some $\sigma=\emptyset, *$ or $P$.

\end{thm}

\noindent{\sl Proof.} Nuzhin~\cite{Nuz90, Nuz97} has shown that a simple group $L_2(q)$ is in $\G(1)$ if and only if $q\ne 7, 9$ (these are the exceptions $L_3(2)$ and $A_6$ in Corollary~\ref{Mcor} and Table~\ref{FSGpsRealised}).  It follows from Lemma~\ref{regmapslemma}(a) that  if $T=2^{\sigma}$, $3$ or $4^{\sigma}$ then $L_2(q)$ is in $\G(T)$ for all $q\ne 7, 9$. We have already shown that $A_6\in\G(T)$ for these classes, so it is sufficient to consider $L_2(7)$.

First let $T=2$. In $L_2(7)\;(=SL_2(7)/\{\pm I\})$ define the involutions
\[s_1=\pm\left(\,\begin{matrix}0&1\cr -1&0\cr \end{matrix}\,\right),
\quad
s_2=\pm\left(\,\begin{matrix}0&2\cr 3&0\cr \end{matrix}\,\right)
\quad{\rm and}\quad
s'=\pm\left(\,\begin{matrix}1&3\cr -3&-1\cr \end{matrix}\,\right).\]
Then $s_1s_2$ and $s_1s'$ have trace $\pm 1$ and hence have order $3$, while $s_2s'$ has trace $\pm 3$ and hence has order~$4$. These involutions generate $L_2(7)$, since no maximal subgroup contains two elements of order $3$ with a product of order~$4$.  Lemma~\ref{class2lemma}(c) now gives $L_2(7)\in\G(T)$ for all $T=2^{\sigma}$, $3$ or $4^{\sigma}$.

If $T=2^{\sigma}{\rm ex}$ or $5^{\sigma}$ then Lemma~\ref{L2qlemma} shows that $L_2(q)\not\in\G(T)$ for all $q$. \hfill$\square$


\section{$L_3(q)$ and $U_3(q)$}

Recall that $L_3(q)$ is simple for all prime powers $q$, and $U_3(q)$ is simple for all prime powers $q>2$. Although $U_3(3)$ is simple, it has special properties which require separate treatment, given later.

\begin{thm}\label{L3qU3qthm}
Suppose that $G=L_3(q)$, or $G=U_3(q)$ with $q>3$. Then $G\in\G(T)$ if and only if either
$T=2^{\sigma}, 3$ or $4^{\sigma}$, or $T=5^{\sigma}$ and $G\ne L_3(2)$.
\end{thm}

\noindent{\sl Proof.} By Corollaries~\ref{Mcor} and~\ref{Lcor}, the simple groups $L_3(q)$ and $U_3(q)$ are not in $\G(T)$ for $T=1$ or $2^{\sigma}{\rm ex}$.

Now let $T=2$, with $N(2)\cong C_2*C_2*C_2$. Malle, Saxl and Weigel~\cite{MSW} have shown that every non-abelian finite simple group $G\not\cong U_3(3)$ is generated by a strongly real element $x$ (one such that $x^a=x^{-1}$ for some involution $a\in G$) and an involution $y$. Thus $G$ is generated by three involutions $a$, $b:=ax$ and $y$. If an automorphism of $G$ fixes the involution $y$ and transposes $a$ and $b$, then it inverts $x$. However, by Theorem~\ref{Lthm} the groups $G=L_3(q)$ and $U_3(q)$ have no generating pairs $x, y$ with this property, so by taking $s_1=a$, $s_2=b$ and $s_3=y$ we see that $G\in\G(T)$ for $T=2$, and hence for $T=2^{\sigma}$, $3$ and $4^{\sigma}$ by Lemma~\ref{class2lemma}(c).

This leaves the case $T=5^{\sigma}$. We will show that if $q>2$ then $L_3(q), U_3(q)\in\G(5^{\sigma})$ for all $\sigma$. For arguments involving Singer cycles and their eigenvalues, we will need Lemma~\ref{priminv}, that in a finite field of order at least $5$, no primitive element is conjugate to its inverse under an automorphism.

First, let $G=L_3(q)$. Then ${\rm Aut}\,G$ is a semidirect product of $P\Gamma L_3(q)$ by a group of order $2$ generated by the graph automorphism (or polarity) $\gamma$ induced by the matrix operation $A\mapsto (A^{-1})^T$. Let $s$ be a Singer cycle in $G$, an element of order $(q^2+q+1)/d$ where $d=(3,q-1)$; its centraliser in $G$ is $\langle s\rangle$, while its centraliser in ${\rm Aut}\,G$ is a cyclic group $C$ of order $q^2+q+1$; this is generated by a Singer cycle in $PGL_3(q)$, which contains $G$ with index $d$. It follows by applying results of Bereczky~\cite{Ber} to $SL_3(q)$ that the only maximal subgroup of $G$ containing $s$ is $M:=N_G(\langle s\rangle)$, a semidirect product of $\langle s\rangle$ by $C_3$.

Simple numerical estimates show that if $q>2$ there exist elements $s'\in G$ such that $s'\not\in M$ (so $\langle s, s'\rangle=G$), $s'$ is not a Singer cycle (so no automorphism can transpose it with $s^{\pm 1}$) and none of the $q^2+q+1$ automorphisms inverting $s$ also inverts $s'$. (It follows from Lemma~\ref{priminv} that the automorphisms inverting $s$ are all conjugate in ${\rm Aut}\,G$ to $\gamma$; the elements of $G$ inverted by $\gamma$ are those corresponding to symmetric or skew-symmetric matrices.) Thus $G\in\G(5^{\sigma})$ for all $\sigma$. However, $L_3(2)\cong L_2(7)\not\in\G(5^{\sigma})$ by Lemma~\ref{L2qlemma}.

A similar argument can be applied to $U_3(q)$. In this case ${\rm Aut}\,G=P\Gamma U_3(q)$, an extension of $PGU_3(q)$, which contains $G$ with index $d=(3, q+1)$, by ${\rm Gal}\,\F_{q^2}$. The analogue in $G$ of a Singer cycle in $L_3(q)$ is an element $s$ of order $(q^2-q+1)/d$; its centraliser in $G$ is $\langle s\rangle$, while its centraliser in ${\rm Aut}\,G$ has order $q^2-q+1$. In each case $s$ is contained in a unique maximal subgroup $M$ of $G$: if $q\ne 3$ or $5$ then $M=N_G(\langle s\rangle)$, a semidirect product of $\langle s\rangle$ by $C_3$, whereas if $q=3$ or $5$ then $M\cong L_2(7)$ or $A_7$. There is a unique conjugacy class $\mathcal A$ of involutions $\alpha\in{\rm Aut}\,G$ inverting such elements $s$, namely the class (denoted by $2B$ for the groups in~\cite{ATLAS}) containing the automorphism $\gamma$ induced by the automorphism $t\mapsto t^q$ of $\F_{q^2}$. Each $s$ is inverted by $q^2-q+1$ automorphisms $\alpha\in\mathcal A$. Each $\alpha\in\mathcal A$ inverts an element $s'\in G$ if and only if $s'\alpha\in\mathcal A$, so the number of elements $s'$ inverted by $\alpha$ is $|{\mathcal A}|=q^2(q-1)(q^2-q+1)$. As before, it follows that there exist such pairs $s, s'$ which generate $G$ with no forbidden automorphisms, so $G\in\G(5^{\sigma})$. \hfill$\square$


\section{$L_4(2^e)$, $U_4(2^e)$, $U_4(3)$ and $U_5(2)$}

Here we deal with the higher-dimensional linear and unitary groups in ${\mathbb L}\cup{\mathbb M}$.

\begin{thm}\label{LU45}
Let $G=L_4(2^e)$, $U_4(2^e)$, $U_4(3)$ or $U_5(2)$. Then $G\in \G(T)$ if and only if $T\ne 1$.
\end{thm}

\noindent{\sl Proof.} By Corollaries~\ref{Mcor} and~\ref{Lcor} these groups are not in $\G(1)$, but they are in $\G(2^*{\rm ex})$ and hence, by Lemma~\ref{regmapslemma}(b), in $\G(T)$ for all $T=2^{\sigma}{\rm ex}, 4^{\sigma}$ and $5^{\sigma}$. The classes $T=2^{\sigma}$ and $3$ remain, and by Lemma~\ref{class2lemma}(c) it is sufficient to show that these groups are in $\G(2)$.

Let $G=L_4(q)$ where $q=2^e$. We may assume that $e>1$, since the group $L_4(2)\cong A_8$ has been dealt with in Theorem~\ref{altgps}. In~\cite[Theorem~2.1]{MSW},  Malle, Saxl and Weigel showed that there are conjugacy classes ${\mathcal C}_1, {\mathcal C}_2\subset G$, containing elements of orders $q^2+1$ and $q^3-1$ respectively, such that the elements of ${\mathcal C}_1$ are strongly real, and for each non-identity conjugacy class ${\mathcal C}_3\subset G$ there are elements $g_i\in{\mathcal C}_i$ ($i=1, 2, 3$) generating $G$ and satisfying $g_1g_2g_3=1$. We can write $g_1=s_2s_3$ for involutions $s_2, s_3\in G$. There are two conjugacy classes of involutions in $G$, of different sizes, depending on the dimensions of the subspaces of the natural module $\F_q^4$ fixed by their elements. The involutions inverting $g_1$, such as $s_2$ and $s_3$, all lie in one of these two classes, so taking ${\mathcal C}_3$ to be the other class of involutions, and putting $s_1=g_3$, we obtain three involutions $s_1, s_2, s_3$ generating $G$, with no automorphism transposing $s_1$ and $s_2$. Thus $G\in\G(2)$, as required.

The same argument can be applied to $U_4(2^e)$, using~\cite[Theorem~2.2]{MSW}, the only difference being that the elements $g_2\in{\mathcal C}_2$ now have order $q^3+1$ rather than $q^3-1$.

For $U_5(2)$, classes ${\mathcal C}_1$ and ${\mathcal C}_2$ in~\cite[Theorem~2.2]{MSW} contain elements of orders $(q^5+1)/d(q+1)=11$ and $q^{4/2}+(-1)^{4/2}=5$ (where $d=(5, q+1)=1$), the latter strongly real. There are two conjugacy classes of involutions, and if we choose ${\mathcal C}_3$ to be the class $2A$ of involutions not inverting elements of order $5$, the proof proceeds as above. 

The group $U_4(3)$, with only one class of involutions and a comparatively large outer automorphism group, resists this line of argument; however, a computer search by Matan Ziv-Av~\cite{Ziv}, using GAP, has shown that it satisfies the conditions of Lemma~\ref{class2lemma}(a), for example with $ab$, $ac$ and $bc$ having orders $4, 5$ and $6$, so Lemma~\ref{class2lemma}(c) also applies here.  \hfill$\square$

\medskip

The last paragraph is the only place where the proof of Theorem~\ref{mainthmsimple} relies on computers. It would be interesting to have a computer-free argument to deal with $U_4(3)$.


\section{$U_3(3)$}

The group $U_3(3)$ was excluded from Theorem~\ref{L3qU3qthm}, so we treat it here as a special case. One could simply carry out a computer search for suitable generating sets, as was done in the preceding section for $U_4(3)$ and $T=2$, but an argument by hand, using properties of this group given in the ATLAS~\cite{ATLAS}, is more satisfying and more instructive.

\begin{thm}\label{U33,43,52}
The group $U_3(3)$ is in $\G(T)$ if and only if $T=3$, $4^{\sigma}$ or $5^{\sigma}$.
\end{thm}

\noindent{\sl Proof.} Let $G=U_3(3)$. Taking $q=3$ in the proof of Theorem~\ref{L3qU3qthm} shows that $G\not\in \G(T)$ for $T=1$ or $2^{\sigma}{\rm ex}$.

We will show that $G\in\G(3)$ by using Lemma~\ref{class3lemma}, which requires a particular generating set of four involutions. There is a single conjugacy class of maximal subgroups $H\cong L_2(7)$ in $G$. Let $s_1$, $s_2$ and $s_3$ be involutions in such a subgroup $H$ corresponding to the elements
\[\pm\left(\,\begin{matrix}-1&1\cr -2&1\cr \end{matrix}\,\right), \quad
\pm\left(\,\begin{matrix}0&1\cr -1&0\cr \end{matrix}\,\right) \quad{\rm and}\quad
\pm\left(\,\begin{matrix}0&3\cr 2&0\cr \end{matrix}\,\right)\]
of $L_2(7)$, so that $s_1s_2$, $s_1s_3$ and $s_2s_3$ have orders $4$, $4$ and $3$ respectively. Since $s:=s_1s_2s_3$ has order $7$, and no proper subgroup of $L_2(7)$ has order divisible by $14$, these involutions generate $H$.

Now $G$ acts as a primitive group of degree $36$ and rank $4$ on the cosets of $H$, with subdegrees $1, 7, 7$ and $21$. The subgroup $\langle s_1, s_2\rangle\cong D_4$ of $H$ thus fixes more than one point, so it lies in a second subgroup $H'\cong L_2(7)$ of $G$. The subgroups of $L_2(7)$ isomorphic to $D_4$ are all conjugate (they are Sylow $2$-subgroups), so as above there is an involution $s_0\in H'$ such that $\langle s_1, s_2, s_0\rangle=H'$ where $s_1s_2$, $s_1s_0$ and $s_2s_0$ have orders $4$, $4$ and $3$. Then $\langle s_0, s_1, s_2, s_3\rangle=\langle H, H'\rangle=G$ and the partitions $02\mid 13$ and $01\mid 23$ satisfy the condition in Lemma~\ref{class3lemma}, so $G\in\G(3)$.

Next, let $T=4$. The involution $s_1$ and the element $s$ of order $7$ defined above generate the maximal subgroup $H$ of $G$. We need to show that there is another involution $s'_2\in G\setminus H$ such that no automorphism of $G$ inverts $s$ and transposes $s_1$ and $s'_2$. There are seven automorphisms inverting $s$, all of them outer automorphisms in the conjugacy class $2B\subset {\rm Aut}\,G=PGU_3(3)$, forming a coset of $C_G(s)=\langle s\rangle$. Now there are $63$ involutions in $G$, of which $21$ are in $H$. We can therefore choose one of the $42$ involutions $s'_2\in G\setminus H$, avoiding the seven images of $s_1$ under the automorphisms inverting $s$. Thus $G\in\G(4^{\sigma})$ for all $\sigma$. 

Finally let $T=5$. There are no maximal subgroups of $G$ containing elements $s$ and $s'$ of orders $7$ and $6$, so any such pair generates $G$. Since they have different orders, it is sufficient to show that there exists such a pair which are not simultaneously inverted by any automorphism. Any element $s$ of order $7$ is inverted by seven automorphisms, all outer, in class $2B$.  Each such automorphism inverts $24$ elements $s'\in G$ of order $6$, so there are at most $168$ such elements $s'$ inverted by automorphisms also inverting $s$. Since $G$ has $504$ elements of order $6$, a suitable pair $s, s'$ exists, so $G\in\G(5^{\sigma})$ for all $\sigma$. \hfill$\square$

\medskip

This complete the proof of Theorem~\ref{mainthmsimple}.


\end{document}